\documentclass{amsart}

\usepackage{amsmath,amssymb,amsthm,amsfonts,enumerate,array, color,lscape,fancyhdr,layout,pst-all}
\usepackage{mathtools}
\usepackage[all,cmtip]{xy}
\usepackage[english]{babel}
\usepackage[latin1]{inputenc}
\usepackage{young}
\usepackage{tikz}
\usepackage{enumitem}
\usepackage{graphicx}
\usepackage{float}
\usepackage{multicol}
\usepackage{ytableau,tikz,varwidth}
\usetikzlibrary{calc}
\usepackage[active]{srcltx}
\usepackage[colorlinks=true,linkcolor=blue,citecolor=blue]{hyperref}
\usepackage{xcolor}
\usepackage{orcidlink}

\newtheorem{theorem}{Theorem}[section]

\newtheorem{lemma}[theorem]{Lemma}
\newtheorem{corollary}[theorem]{Corollary}
\newtheorem{prop}[theorem]{Proposition}

\newtheorem{defin}[theorem]{Definition}

\newtheorem{teo}[theorem]{Theorem}

\newlength\mylen
\newlist{mycases}{enumerate}{1}
\setlist[mycases,1]{label=\textit{Case~\arabic*.}, 
  labelwidth=\dimexpr-\mylen-\labelsep\relax,leftmargin=0pt,align=right}

\makeatletter
\DeclarePairedDelimiterX{\pmodx}[1]{(}{)}{{\operator@font mod}\mkern6mu#1}
\renewcommand{\pmod}{%
  \allowbreak
  \if@display\mkern18mu\else\mkern8mu\fi
  \pmodx
}
\makeatother

\begin{document}
\title[Minimal product set in non-abelian metacyclic groups of even order]{Minimal product set in non-abelian metacyclic groups of even order}

\author[F. A. Benavides]{Fernando A. Benavides \orcidlink{0000-0002-4115-8408}}
\address{Departamento de Matem\'aticas y Estad\'istica \\ Universidad de Nari\~no \\ Pasto \\ Colombia\\}
\email{fandresbenavides@udenar.edu.co}

\author[W. Mutis]{Wilson F. Mutis \orcidlink{0009-0004-2664-4362}}
\address{Departamento de Matem\'aticas y Estad\'istica \\ Universidad de Nari\~no \\ Pasto \\ Colombia}
\email{wilsonmutis@udenar.edu.co}

\begin{abstract}
Given a finite group $G$ and positive integers $r$ and $s$, a problem of interest in algebra is determining the minimum cardinality of the product set $AB$, where $A$ and $B$ are subsets of $G$ such that $|A|=r$ and $|B|=s$. This problem has been solved for the class of abelian groups; however, it remains open for finite non-abelian groups. In this paper, we prove that the result obtained for abelian groups can be extended to the class of metacyclic groups $K_{m,n}=\left\langle a,b \ : \ a^m=1,b^{2n}=a^g,bab^{-1}=a^{-1}\right\rangle$. Consequently, we provide a new proof of the result for the dihedral group $D_n$ and dicylic group $Q_{4n}$.
\end{abstract}


\maketitle

\section{Introduction}

Let $G$ be a group and $A,B$ be subsets of $G$. Minimizing the size of product set $AB$ under the conditions $|A|=r$ and $|B|=s$ for arbitrary integers $r,s$ with $1\leq r,s\leq |G|$ is a computationally difficult problem. Therefore, an interesting problem in additive number theory and group theory is to model
\begin{equation*}
\mu_G\left(r,s\right)=\min\left\{\left\vert AB\right\vert:A,B\subseteq G,\left\vert A\right\vert=r,\left\vert B\right\vert=s\right\}.
\end{equation*}
using a simpler function. Cauchy \cite{cauchy1812} and Davenport independently \cite{davenport1935} studied $\mu_G$ for the cyclic group of prime order, stating that
\begin{equation*}
\mu_{\mathbb{Z}_p}(r,s)=\min\{r+s-1,p\},
\end{equation*}
with $1\leq r,s\leq p$. Eliahou and Kervaire \cite{Eliahou2005} established a relationship between $\mu_G$ and the arithmetic function $\kappa_G$, given by
\begin{equation*}
\kappa_G\left(r,s\right)=\min_{h\in \mathcal{H}(G)} h\left( \left\lceil  {\frac{r}{h}} \right \rceil+\left\lceil{\frac{s}{h}} \right \rceil-1\right),
\end{equation*}
where $\mathcal{H}(G)$ denotes the set of the subgroup orders of $G$. They proved that $\mu_G=\kappa_G$ holds true if $G$ is a finite abelian group. Eliahou and Kervaire \cite{Eliahou2010}, Mutis, Benavides and Castillo \cite{Fernando2010,Mutis2012},  Berchenko-Kogan \cite{Berko}  and recently Kaur and Singh \cite{kaur-sing2024} proved that equality $\mu_G=\kappa_G$ also is satisfied if $G$ is the dihedral group $D_n$ of order $2n$, a finite $p$-group, a finite hamiltonian group or dicyclic group $Q_{4n}$ of order $4n$.

\vspace{0.3cm}

Although there are similar results for particular cases in non-abelian groups, the general case remains open.  However, Eliahou and Kervaire \cite{Eliahou2010} showed that in the semi-direct product $G=C_7\rtimes C_3$ the equality $\mu_G(r,s)=\kappa_G(r,s)$ fails if $(r,s)\in\{(5,9),(6,8),(6, 9),(8, 9),(9, 9)\}$. Using a modified version of the arithmetic function $\kappa_G$, Deckelbaum \cite{Alan} established the values of $\mu_G$ for non-abelian groups $G$ of the order $3p$. This result was later extended by Berchenko-Kogan \cite{Berko} to all non-abelian groups of order $pq$ with $p>q$ odd primes. He proved that
\begin{equation*}
\mu_G(r,s) = \left\lbrace
\begin{array}{ll}
\mathcal{N}\kappa_G(r,s)  & \textup{if } \ r,s>q \ \text{and} \ \left\lceil \frac{r}{q}\right \rceil + \left\lceil \frac{s}{q}\right\rceil < p \\\\
\kappa_G(r,s) & \textup{other case.} 
\end{array}
\right.
\end{equation*}
Function $\mathcal{N}\kappa_G$ was introduced by Eliahou and Kervaire \cite{Eliahou2010}.

\vspace{0.3cm}

In this paper, we study the function $\mu_G$ for non-abelian groups $K_{m,n}$  where 
\begin{equation*}
K_{m,n}=\left\langle a,b \ : \ a^m=1,b^{2n}=a^g,bab^{-1}=a^{-1}\right\rangle,
\end{equation*}
and $m,n\geq 1$. This group was studied by Kaur and Singh \cite{kaur-sing2024} in case that $g=0$. We characterize the divisors $h$ of $2mn$ for which there is no a normal subgroup of order $h$ and prove by double induction that $\mu_G(r,s)=\kappa_G(r,s)$ if $G=K_{m,n}$ and $1\leq r,s,\leq 2mn$. Consequently, $\mu_G=\kappa_G$ holds true if $G=D_n$ or $G=Q_{4n}$.

\section{Preliminaries} \label{sec:prel}

In this section, we give some definitions, notations, and results needed for our main proof. In the following, $G$ always refers to a finite multiplicative group.
\vspace{0.3cm}

\subsection{Supporting results and notations}

Given two non-empty subsets $A$ and $B$, the product $AB$ is the set of elements of $G$ which can be written in at least one way as a product of an element of $A$ with an element of $B$, that is 
\begin{equation*}
AB=\left\{ab:a\in A, b\in B\right\}.
\end{equation*} 
\begin{defin}
Let $G$ be a finite group and let $r$ and $s$ be positive integers such that $1\leq r,s\leq |G|$. The function $\mu_G(r,s)$ is the minimum cardinality of the product of a subset of $G$ of cardinality $r$ and a subset of $G$ of cardinality $s$. That is,
\begin{equation*}
\mu_G\left(r,s\right)=\min\left\{\left\vert AB\right\vert:A,B\subseteq G,\left\vert A\right\vert=r,\left\vert B\right\vert=s\right\}.
\end{equation*}
\end{defin}
The minimal product set problem aims to characterize the values of the function $\mu_G(r,s)$. If $A$ and $B$ are subsets of $G$ such that $\mu_G(r,s)=|AB|$ with $|A|=r$ and $|B|=s$, then we say that $A$ and $B$ realize $\mu_G(r,s)$.

\vspace{0.5cm}

The function $\mu_G$ has been studied for a few classes of groups. For example, in \cite{Eliahou2010} the authors bounded $\mu_G$ by means the functions
\begin{center}
$\mathcal{D}\kappa_G(r, s)=\min\limits_{h\in\mathcal{D}(G)} f_h(r,s)$ \ and \ $\mathcal{N}\kappa_G(r, s)=\min\limits_{h\in\mathcal{N}(G)} f_h(r,s)$,
\end{center}
where $\mathcal{D}(G)$ is the set of positive divisors of $|G|$, $\mathcal{N}(G)$ is the set of orders of normal subgroups of $G$ and
\begin{equation*}
f_h(r,s)= \lceil r\rceil_h+\lceil s\rceil_h-h,
\end{equation*}
where $\lceil r\rceil_h=\min\left\{n\in\mathbb{Z}^+:r\leq n, \ n\equiv 0\hspace{-0.1cm}\pmod{h}\right\}$.

\begin{teo}[Eliahou and Kervaire \cite{Eliahou2010}]\label{teor-soluble}
Let $G$ be a finite solvable group and $1\leq r,s\leq\left\vert G\right\vert$. Then 
\begin{equation*}
\mathcal{D}\kappa_G(r, s)\leq\mu_G(r, s)\leq\mathcal{N}\kappa_G(r, s).
\end{equation*}
\end{teo}
Note that if a finite group $G$ satisfies the converse of Lagrange's Theorem, then
\begin{equation}\label{soluble-reciproco}
\kappa_G(r, s)\leq\mu_G(r, s)\leq\mathcal{N}_{\kappa_G}(r, s).
\end{equation}

The Dihedral groups $D_n$ of order $2n$, finite $p$-groups, finite Hamiltonian groups, and non-abelian groups of order $pq$, where $p$ and $q$ are distinct odd primes, are finite solvable groups to which the relation (\ref{soluble-reciproco}) applies. This argument plays a fundamental role in the proof of the main results of Eliahou and Kervaire \cite{Eliahou2006}, Mutis, Benavides, and Castillo \cite{Fernando2010,Mutis2012},  Deckelbaum \cite{Alan}, and Berchenko-Kogan \cite{Berko}.

\subsection{Metacyclic groups}

A finite group $K$ is metacyclic if there exists a normal subgroup $A\unlhd K$ such that $A$ and the quotient group $K/A$ are both cyclic groups. According to Yang \cite{Yang2020}, we assume that all finite metacyclic groups $K$ admit the following presentation
\begin{equation*}
K=\left\langle a,b:a^m=1,b^n=a^g, bab^{-1}=a^h\right\rangle
\end{equation*}
where $m$, $n$, $g$ and $h$ are integers that satisfy the following relations
\begin{enumerate}
\item $m\geq1$ and $n\geq1$.
\item $0\leq g,h\leq m-1$.
\item $g(h-1)\equiv 0\pmod m$ and $h^n\equiv1\pmod m$.
\end{enumerate}
Consider the set $\Omega$ of all $4$-tuples $\left(m,n,g,h\right)\in\mathbb{N}^4$ which verify the three aforementioned relations. We present the version of H\"older theorem given in \cite{Yang2020}, which will be used throughout this document.

\begin{teo}[H\"older]\label{holder}
Let $G$ be a finite metacyclic group with cyclic normal subgroup $A$ and cyclic quotient $G/A$. Then there exists $\left(m,n,g,h\right)\in\Omega$, such that $|A|=m$, $|G/A|=n$, and
\begin{equation*}
G\cong\left\langle a,b:a^m=1,b^n=a^g, bab^{-1}=a^h\right\rangle.
\end{equation*}
Furthermore, fix a $\left(m,n,g,h\right)\in\Omega$ and let  $K=\left\langle a,b:a^m=1,b^n=a^g, bab^{-1}=a^h\right\rangle.$ Then we have the following valid statements.
\begin{enumerate}
\item $\left\vert a\right\vert=m$ and $\langle a\rangle$ is a normal subgroup of $K$.
\item $K/\left\langle a\right\rangle=\left\langle \overline{b}\right\rangle$, $\left\vert \overline{b}\right\vert=n$, and $K$ is a metacyclic group of order $mn$.
\item $\langle a\rangle\cap\left\{b^j: j=1,\ldots,n-1\right\}=\emptyset.$
\item For all pairs $(i,j),(r,s)\in\{0,1,\ldots,m-1\}\times\{0,1,\ldots,n-1\}$, we have
$$a^ib^j=a^rb^s\text{, if and only if, }(i,j)=(r,s).$$
\end{enumerate}
\end{teo}

Metacyclic groups satisfy the converse of Lagrange's  theorem and are solvable, implying that,
\begin{equation*}
\kappa_K(r, s)\leq\mu_K(r, s)\leq\mathcal{N}\kappa_K(r, s).
\end{equation*}
To study the function $\mu_K$, it is necessary to determine the set $\mathcal{N}(K)$ of orders of the normal subgroups. Theorem 3.2 of \cite{Yang2020}, Yang states a bijection between the sets
\begin{equation*}
\Gamma=\left\{(k,l,\beta)\in\mathbb{N}^3 \ : \ l\mid m, \ l\mid k, \ k\mid nl, \ \beta<\dfrac{m}{l}, \ \beta\sum\limits_{j=0}^{\frac{k}{l}-1}h^{\frac{nl}{k}j}\equiv-g\hspace{-0.3cm}\pmod[\Big]{\frac{m}{l}} \right\}.
\end{equation*}
and $T=\left\{A: A\leq K\right\}$, given by
\begin{equation*}
\Psi(k,l,\beta)=\left\langle a^{\frac{m}{l}},a^{\beta}b^{\frac{nl}{k}}\right\rangle.
\end{equation*}

Important properties that $\Psi$ satisfies are
\begin{center}
$|\Psi(k,l,\beta)|=k$ \ and \ $|\Psi(k,l,\beta)\cap\langle a\rangle|=l$. 
\end{center}
Further, in Theorem 3.3, he provides a necessary and sufficient condition to determine whether a subgroup $A$ is normal in $K$. These results are summarized in the following theorem:
\begin{teo}[Yang \cite{Yang2020}]\label{yang-teor}
Let $K$ be the metacyclic group
\begin{equation*}
K=\left\langle a,b:a^m=1,b^n=a^g, bab^{-1}=a^h\right\rangle,
\end{equation*}
where $(m,n,g,h)\in\Omega$. Then,
\begin{enumerate}
\item The map $\Psi:\Gamma\longrightarrow T$ given by
\begin{equation*}
\Psi(k,l,\beta)=\left\langle a^{\frac{m}{l}},a^{\beta}b^{\frac{nl}{k}}\right\rangle,
\end{equation*}
is a bijection.
\item $\Psi\left(\rho,l,\beta\right)\unlhd K$, if and only if,
\begin{equation*}
\beta(h-1)\equiv 0\hspace{-0.3cm}\pmod[\Big]{\frac{m}{l}},\text{ and }h^{\frac{nl}{\rho}}\equiv1\hspace{-0.3cm}\pmod[\Big]{\frac{m}{l}}.
\end{equation*} 
\item If $\Psi\left(\rho,l,\beta\right)\unlhd K$, then $\left(r=\dfrac{m}{l},s=\dfrac{nl}{\rho},t=-\beta\%\dfrac{m}{l},u=h\%\dfrac{m}{l}\right)\in\Omega$, with $p\%q$ denote the only integer $w\in\{0,1,\ldots,\vert q\vert-1\}$ such that $p\equiv w \pmod q$. Furthermore, 
\begin{equation*}
K/\Psi\left(\rho,l,\beta\right)\cong\left\langle x,y:x^r=1,y^s=x^t, yxy^{-1}=x^u\right\rangle
\end{equation*} 
\end{enumerate}
\end{teo}

\section{Main theorem}

For positive integers $m$ and $n$, consider the metacyclic group
\begin{equation*}
K_{m,n}=\left\langle a,b \ : \ a^m=1,b^{2n}=a^g,bab^{-1}=a^{-1}\right\rangle.
\end{equation*}
where $(m,n,g,m-1)\in\Omega$. Observe that, if $m$ is odd, then $g=0$; otherwise, $g=0$ or $\frac{m}{2}$. From now on, $K_{m,n}$ denotes a non-abelian group.

\begin{lemma}\label{lem:normal-impar}
Let $m,n\geq 1$ be integers. If $k$ is an odd positive divisor of $2mn$, then there is a $(k,l,\beta)\in\Gamma$ such that $\Psi(k,l,\beta)\trianglelefteq K_{m,n}$. 
\end{lemma}
\begin{proof}
Let $(k,l,\beta)$ be an element of $\Gamma$ such that $k\equiv 1\pmod{2}$. Then
\begin{equation*}
v_2(2nl)-v_2(k)>0,
\end{equation*}
where $v_2(\cdot)$ is the $2$-adic valuation. If $g=0$ then $(k,l,0)\in\Gamma$ and $\Psi(k,l,0)\trianglelefteq K_{m,n}$. In case that $m$ is even and $g=\frac{m}{2}$, $\left(k,l,\frac{m}{2l}\right)\in\Gamma$ since
\begin{equation*}
\left(\dfrac{m}{2l}\right)\left(\dfrac{k}{l}\right)+\dfrac{m}{2}=\left(\dfrac{m}{2l}\right)\left(\dfrac{k}{l}+l\right)\equiv 0\hspace{-0.3cm}\pmod[\Big]{\dfrac{m}{l}},
\end{equation*}
and besides, $\Psi\left(k,l,\frac{m}{2l}\right)\trianglelefteq K_{m,n}$.
\end{proof}

\begin{teo}\label{teor-v2}
Let $k$ be a positive divisor of $2mn$. 
\begin{center}
$k\in\mathcal{N}(K_{m,n})$, if and only if, $m\mid k$ or $v_2(k)<v_2(2mn)$.
\end{center}
\end{teo}
\begin{proof}
According to Lemma~\ref{lem:normal-impar}, we can assume that $k$ is even and $m\nmid k$. Then, there exists $(k,l,\beta)\in\Gamma$ such that $\Psi(k,l,\beta)\trianglelefteq K_{m,n}$. If $v_2(k)=v_2(2nl)$ then $m=2l$, and so
\begin{equation*}
v_2(2mn)>v_2(mn)=v_2(k).
\end{equation*}

Conversely, if $m\mid k$, then $(k,m,0)\in\Gamma$ and $\Psi(k,m,0)\trianglelefteq K_{m,n}$. Let $(k,l,\beta)\in\Gamma$ be a tuple, such that $m\nmid k$. If $v_2(2nl)-v_2(k)>0$, then $(k,d,0)\in\Gamma$ and $\Psi(k,d,0)\trianglelefteq K_{m,n}$ where $d=\gcd(k,m)$. However, if $v_2(k)=v_2(2nl)$ then $v_2(l)<v_2(m)$, implying that $2\mid m$. Therefore, $(k,2l,0)\in\Gamma$ and $\Psi(k,2l,0)\trianglelefteq K_{m,n}$.
\end{proof}

\begin{prop}\label{prop:kappa-mu-f2}
Let $m,n\geq 1$ be odd integers and $1\leq r,s\leq 2mn$. If $\kappa_G(r,s)=f_2(r,s)$ with $G=K_{m,n}$, then 
\begin{equation*}
\mu_G(r,s)=\kappa_G(r,s),
\end{equation*}
\end{prop}
\begin{proof}
Since $f_2(r,s)<f_1(r,s)$, $r$ and $s$ must be even integers. Suppose that $r=2r_1$ and $s=2s_1$ then, we apply the Euclidean division on $r_1$ and $s_1$ by $m$. We have  $r_1=mp+x$ and $s_1=mq+y$ with $0\leq x,y<m$. Note that in the case $x=y=0$, $f_m(r,s)< f_2(r,s)$, which is not possible.

\vspace{0.5cm}

Let $F=\langle a\rangle$ be the normal subgroup of $G$ generated by $a$. For $0\leq i<n$, consider the subset $F_i=b^{2i}F=Fb^{2i}$, and for each positive integer $q$, the initial segment of cosets $I(q)$ is defined by
\begin{equation*}
I(q)=F_0\cup F_1\cup\cdots\cup F_{q-1}.
\end{equation*}
Under this notation we have $I(p)I(q)=I(p+q-1)$. For $x\geq 1$ define $K_x=\{1,a,\ldots,a^{x-1}\}$ and empty in other case. Consider the following cases:
\begin{mycases}
\item ($p=q$) If $A'=I(p)b^2\cup K_x$ and $B'=I(p)b^2\cup K_y$ then,
\begin{equation*}
A'B'\subseteq I(2p)b^2\cup K_xK_y.
\end{equation*}

\item ($p<q$) If $A'=I(p)\cup b^{2p}K_x$ and $B'=I(p)\cup K_yb^{2q}$ then,
\begin{equation*}
A'B'\subseteq I(p+q)b^2\cup K_xK_yb^{2(p+q)}.
\end{equation*}
\end{mycases}
In both cases, we have $|A'|=r_1$, $|B'|=s_1$ and $|A'B'|\leq r_1+s_1-1$. This inequality holds even if $x=0$ and $y\geq 1$, in which case $K_x=\emptyset$ and $K_xK_y=\emptyset$.

\vspace{0.5cm}

Finally, for the subsets $A=A'\cup A'a^{y-1}b^n$ and $B=B'\cup B'b^n$, it is satisfied that $AB=A'B'\cup A'B'b^n$ and, in consequence,
\begin{equation*}
|AB|\leq 2|A'B'|\leq 2(r_1+s_1-1)=r+s-2=f_2(r,s).
\end{equation*}
\end{proof}

\begin{teo}
Let $m,n\geq 1$ be integers. If $G=K_{m,n}$ then 
\begin{equation*}
\mu_{G}(r,s)=\kappa_{G}(r,s)
\end{equation*}
for all $1\leq r,s\leq2mn$.
\end{teo}
\begin{proof}
Let $h\in\mathcal{H}(G)$ be the order of a subgroup of $G$ such that $\kappa_{G}(r,s)=f_h(r,s)$. If there is a normal subgroup $H$ of $G$ such that $|H|=h$ or $h=2$, then, by Theorem~\ref{teor-soluble} or Proposition~\ref{prop:kappa-mu-f2}, we have $\mu_{G}(r,s)=\kappa_{G}(r,s)$. 

\vspace{0.5cm}

Otherwise, Theorem~\ref{teor-v2} implies that $h=2k$ with $k>1$, $m\nmid 2k$ and $v_2(k)=v_2(mn)$. In this case, the proof is done by double induction on $m+n$. Let $(k,l,\beta)\in\Gamma$ be a tuple such that $\Psi(k,l,\beta)=K\trianglelefteq G$, then by Theorem~\ref{yang-teor}, the quotient group $U_{m,n}/K$ has the following presentation:
\begin{equation*}
U_{m,n}/K\cong\left\langle c,d \ : \ c^{\overline{m}}=1, \ d^{2\overline{n}}=c^{\overline{g}}, \ dcd^{-1}=c^{-1}\right\rangle=K_{\overline{m},\overline{n}},
\end{equation*}
where $\overline{m}=\dfrac{m}{l}$ and $\overline{n}=\dfrac{nl}{k}$. Note that $\overline{m}+\overline{n}<m+n$ then, by the induction hypothesis
\begin{equation*}
\kappa_{U_{\overline{m},\overline{n}}}(2z,2w)=\mu_{U_{\overline{m},\overline{n}}}(2z,2w),
\end{equation*}
where $w=\left\lceil \dfrac{r}{2k}\right\rceil$ and $z=\left\lceil \dfrac{s}{2k}\right\rceil$. Let $\overline{S}$ and $\overline{T}$ be subsets of $K_{\overline{m},\overline{n}}$ such that $|\overline{S}|=2w$, $|\overline{T}|=2z$, and
\begin{equation*}
|\overline{S}\overline{T}|\leq f_2(2w,2z).
\end{equation*}
Therefore, if $S=K\overline{S}$ and $T=K\overline{T}$ then, both are subsets of $K_{m,n}$ and 
\begin{equation*}
|ST|=k|\overline{S}\overline{T}|\leq k(2w+2z-2)=f_{2k}(2w,2z).
\end{equation*}
Because $f_{2k}(r,s)<f_k(r,s)$, we have $\lceil r\rceil_{2k}=\lceil r\rceil_{k}$ and $\lceil s\rceil_{2k}=\lceil s\rceil_{k}$. This implies that $r\leq 2wk$ and $s\leq 2zk$. In consequence, there exists subsets $A\subseteq S$ and $B\subseteq T$ such that $|A|=r$, $|B|=s$ and
\begin{equation*}
|AB|\leq f_{h}(r,s)=\kappa_G(r,s).
\end{equation*}
\end{proof}

\begin{corollary}
Let $n$ be a positive integer. Then,
\begin{enumerate}
\item $\mu_{Q_{4n}}(r,s)=\kappa_{Q_{4n}}(r,s)$ if $Q_{4n}$ is a dicyclic group of order $4n$ and $1\leq r,s\leq 4n$.
\item $\mu_{D_{2n}}(r,s)=\kappa_{D_{2n}}(r,s)$ if $D_{4n}$ is a dihedral group of order $2n$ and $1\leq r,s\leq 2n$.
\end{enumerate}
\end{corollary}

\section{Acknowledgments}

This work was supported by Vicerector\'ia de Investigaci\'on e Interacci\'on Social of Universidad de Nari\~no. 


\bibliographystyle{plain}
\bibliography{bibliografia}

\end{document}